# DEFINITIONS OF COMPLEX ORDER INTEGRALS AND COMPLEX ORDER DERIVATIVES USING OPERATOR APPROACH

[1] **Raoelina Andriambololona,** [2] **Ranaivoson Tokiniaina,** [3] **Hanitriarivo Rakotoson**
[1] Theoretical Physics Department, Institut National des Sciences et Techniques Nucléaires (INSTN-Madagascar)
instn@moov.mg, raoelinasp@yahoo.fr, jacquelineraoelina@yahoo.fr
[2] INSTN-Madagascar, tokhiniaina@gmail.com
[3] INSTN-Madagascar, infotsara@gmail.com

*Abstract*— For a complex number $s$ , the $s$-order integral of a function $f$ fulfilling some conditions is defined as the action of an operator, noted $J^s$, on $f$. The definition of the operator $J^s$ is given firstly for the case of complex number $s$ with positive real part. Then, using the fact that the operator of one order derivative, noted $D^1$, is the left hand side inverse of the operator $J^1$, an $s$-order derivative operator, noted $D^s$, is also defined for complex number $s$ with positive real part. Finally, considering the relation $J^s = D^{-s}$, the definition of the $s$-order integral and $s$-order derivative is extended for any complex number $s$. An extension of the definition domain of the operators is given too.

## I. INTRODUCTION

The problem of the extension of real integer-order derivatives to fractional derivatives is an old one introduced by Leibnitz in 1695. Several approaches have been done and bibliography may be found in references [1], [2]

In our work [3], we have given an unified definition for both derivatives and indefinite integrals of any power function $f(x) = ax^k$ defined on $\mathbb{R}_+$. Notions of linear, semi-linear, commutative and semi-commutative properties of fractional derivatives have been introduced too.

Our investigations are going on. In another work [4], we have utilized a more general approach. We have considered a larger set of causal functions, which is the definition domain of $s$-order integral operators $J^s$. The $k$-order derivative operator $D^k$ is derived from the $s$-order integral operators $J^s$. Properties of $J^s$ and $D^k$ for any positive real $s$ and positive real $k$ have been studied. Remarkable relations verified by $J^\pi, J^e, D^\pi$ and $D^e$ for transcendental numbers $\pi$ and $e$ have been given too.

We have studied the case of $s$ and $k$ real numbers (case of $\mathbb{R}$ field)[4].

In the present paper, we will give the extension to the case of complex numbers ($\mathbb{C}$ field). The extension is neither trivial nor straightforward because we have to study the existence of the $s$ -order integral operator $J^s$ and the $k$ -order derivative operator $D^k$ in the case of complex numbers $s$ and $k$.

It was shown in [4] that for a derivable and integrable function $f$ of real variable verifying the condition

$$f(x) = 0 \text{ for } x \le 0 \qquad (1.1)$$

and for $s \in \mathbb{R}_+$, a $s$-order integral and $s$-order derivative of the function $f$ can be defined respectively as

$$J^s(f)(x) = \frac{1}{\Gamma(s)} \int_0^x (x-y)^{s-1} f(y) dy \qquad (1.2a)$$

$$= \frac{x^s}{\Gamma(s)} \int_0^1 (1-u)^{s-1} f(ux) du \qquad (1.2b)$$

$$D^s(f)(x) = D^k J^{k-s}(f)(x) \qquad (1.3)$$

with $k > s$. $\Gamma$ is the Euler gamma function [5], [6]. Now we extend those results for $s \in \mathbb{C}$.

*Remark*
We would like to point out that throughout our work, we utilize the notation for the operator product i.e

$$J^a(J^b)(f)(x) = J^a[J^b[f(x)]]$$

$B(p)(q)$ stands for the writing of the usual $B(p, q)$

## II. DEFINITION OF $s$-ORDER INTEGRAL OPERATOR $J^s$ FOR $s$ A COMPLEX NUMBER WITH $Re(s) > 0$

*Theorem 1*

For a function $f$ verifying the conditions $f(x) = 0$ for $x \le 0$, the integral $J^s(f)(x)$ defined by the relation

$$J^s(f)(x) = \frac{1}{\Gamma(s)} \int_0^x (x-y)^{s-1} f(y) dy \qquad (2.1)$$

is absolutely convergent for $s \in \mathbb{C}$ and $Re(s) > 0$.

*Proof*

Let $s = \alpha + i\beta$ then

$$\left| \frac{(x-y)^{s-1}}{\Gamma(s)} f(y) \right| = \frac{|x-y|^{\alpha-1}}{|\Gamma(s)|} |f(y)|$$

$$\Rightarrow \int_0^x \left| \frac{(x-y)^{\alpha-1}}{\Gamma(s)} f(y) \right| dy = \int_0^x \frac{(x-y)^{\alpha-1}}{|\Gamma(s)|} |f(y)| dy$$

Using an integration by parts, we have for $Re(s) = \alpha > 0$

$$\int_0^x \frac{(x-y)^{\alpha-1}}{|\Gamma(s)|} |f(y)| dy$$

$$= -\frac{1}{\alpha|\Gamma(s)|} \{ \lim_{y \to x^-} (x-y)^{\alpha} |f(y)| - \lim_{y \to 0^+} (x-y)^{\alpha} |f(y)| $$

$$- \int_0^x (x-y)^{\alpha} |f(y)|' dy \}$$

$$\lim_{y \to x^-} (x-y)^{\alpha} |f(y)| = 0 \qquad \lim_{y \to 0^+} (x-y)^{\alpha} |f(y)| = 0$$

So

$$\int_0^x \frac{(x-y)^{\alpha-1}}{|\Gamma(s)|} |f(y)| dy = \int_0^x \frac{(x-y)^{\alpha}}{\alpha|\Gamma(s)|} |f(y)|' dy$$

For $x \in I = [0, +\infty[$ and for $\alpha > 0$

$$0 \leq y \leq x \Rightarrow 0 \leq x - y \leq x \Rightarrow 0 \leq (x-y)^{\alpha} \leq x^{\alpha}$$

$$\Rightarrow 0 \leq \frac{(x-y)^{\alpha}}{\alpha|\Gamma(s)|} \leq \frac{x^{\alpha}}{\alpha|\Gamma(s)|}$$

$$\Rightarrow \int_0^x \frac{(x-y)^{\alpha}}{\alpha|\Gamma(s)|} |f(y)|' dy \leq \frac{x^{\alpha}}{\alpha|\Gamma(s)|} \int_0^x |f(y)|' dy$$

$$\Rightarrow \int_0^x \frac{(x-y)^{\alpha}}{\alpha|\Gamma(s)|} |f(y)|' dy \leq \frac{x^{\alpha}}{\alpha|\Gamma(s)|} |f(x)|$$

$$\Rightarrow \int_0^x \left| \frac{(x-y)^{s-1}}{\Gamma(s)} f(y) \right| dy$$

$$= \int_0^x \frac{(x-y)^{\alpha}}{\alpha|\Gamma(s)|} |f(y)|' dy \leq \frac{x^{\alpha}}{\alpha|\Gamma(s)|} |f(x)|$$

Then for $\alpha = Re(s) > 0$, we obtain

$$\int_0^x \left| \frac{(x-y)^{s-1}}{\Gamma(s)} f(y) \right| dy \leq \frac{x^{\alpha}}{\alpha|\Gamma(s)|} |f(x)| \quad (2.2)$$

The integral $J^s(f)(x)$ is absolutely convergent for $\alpha =$

$Re(s) > 0$.

As an extension of the relation(1.2) we consider the following definition.

### Definition 1

Let $f$ be a function which fulfills the condition $f(x) = 0$ for $x \leq 0$. For complex number $s$ with $Re(s) > 0$, we define the $s$-order integral of the function $f$, considered as the action of the operator $J^s$ on this function, by the relation

$$J^s(f)(x) = \frac{1}{\Gamma(s)} \int_0^x (x-y)^{s-1} f(y) dy \quad (2.3)$$

$$= \frac{x^s}{\Gamma(s)} \int_0^1 (1-u)^{s-1} f(ux) du$$

according to the theorem 1, the integral $J^s(f)(x)$ is well defined.

*Example*

$$f(x) = x^{a+ib} = e^{(a+ib)\ln(x)} \qquad a + ib \in \mathbb{C}$$

$$= x^a [\cos(b\ln x) + i\sin(b\ln x)]$$

$$J^s(f)(x) = \frac{1}{\Gamma(s)} \int_0^x (x-y)^{s-1} f(y) dy$$

$$= \frac{x^s}{\Gamma(s)} \int_0^1 (1-u)^{s-1} f(ux) du$$

$$J^s(x^{a+ib}) = \frac{x^{s+a+ib}}{\Gamma(s)} \int_0^1 (1-u)^{s-1} u^{a+ib} du$$

$$= \frac{x^{s+a+ib}}{\Gamma(s)} B(s, a + ib + 1)$$

in which $B$ is the extension of the Euler's beta function for complex numbers

$$B(s, a + ib + 1) = \int_0^1 (1-u)^{s-1} u^{a+ib} du$$

$$= \frac{\Gamma(s)\Gamma(a + ib + 1)}{\Gamma(s + a + ib + 1)}$$

in which $\Gamma$ is the extension of Euler's gamma function for complex numbers

$$J^s(x^{a+ib}) = \frac{x^{s+a+ib}}{\Gamma(s)} \frac{\Gamma(s)\Gamma(a + ib + 1)}{\Gamma(s + a + ib + 1)}$$

$$= \frac{\Gamma(\alpha + 1 + i\beta)}{\Gamma(s + a + ib + 1)} x^{s+a+ib} \quad (2.4)$$

If $s = \alpha + i\beta$ in the relation (2.4), we have

$$J^s(x^{a+ib}) = \frac{\Gamma(a+ib+1)}{\Gamma(\alpha+i\beta+a+ib+1)}x^{a+\alpha+ib+i\beta}$$

$$= \frac{\Gamma(a+ib+1)}{\Gamma(\alpha+a+1+i\beta+ib)}x^{a+\alpha}$$
$$. \left[cos((\beta+b)lnx) + isin((\beta+b)lnx)\right] \quad (2.5)$$

From the relation (2.4), it can be easily deduced that

$$J^{s_1}J^{s_2}(x^{a+ib}) = J^{s_1+s_2}(x^{a+ib}) = J^{s_2}J^{s_1}(x^{a+ib}) \quad (2.6)$$

## III. Properties of the operator $J^s$ for $Re(s) > 0$

### Theorem 2

The set $E$ of function $f$ verifying the conditions $f(x) = 0$ if $x \leq 0$ is a vectorial space on the field $\mathbb{C}$ of complex number and the operator $J^s$ is a linear operator on $E$.

### Proof
The set $E$ of functions which verify the conditions (1.1) is a vectorial space on the field $\mathbb{C}$:

$$\forall (f,g) \in E^2 : f(x) = 0 \text{ and } g(x) = 0 \text{ if } x \leq 0$$
$$\Rightarrow \forall(\lambda,\mu) \in \mathbb{C}^2 : \lambda f(x) + \mu g(x) = 0 \text{ if } x \leq 0$$
$$\Rightarrow (\lambda f + \mu g) \in E$$

The operator $J^s$ is a linear operator on $E$

$$\forall f \in E : J^s(f)(x) = \frac{1}{\Gamma(s)}\int_0^x (x-y)^{s-1}f(y)dy$$

$$f(x) = 0 \text{ if } x \leq 0 \Rightarrow J^s(f)(x) = 0 \text{ if } x \leq 0$$
$$\Rightarrow J^s(f) \in E$$

$$\Rightarrow \forall(\lambda,\mu) \in \mathbb{C}^2, \ \forall (f,g) \in E^2 :$$

$$J^s(\lambda f + \mu g)(x) = \frac{1}{\Gamma(s)}\int_0^x (x-y)^{s-1}[\lambda f(y) + \mu g(y)]dy$$

$$= \frac{\lambda}{\Gamma(s)}\int_0^x (x-y)^{s-1}f(y)dy + \frac{\mu}{\Gamma(s)}\int_0^x (x-y)^{s-1}g(y)dy$$

$$= \lambda J^s(f)(x) + \mu J^s(g)(x) \quad (3.1)$$

### Theorem 3

For all complex numbers $s_1$ and $s_2$ with $Re(s_1) > 0$ and $Re(s_2) > 0$, we have the semi-group property

$$J^{s_1}J^{s_2}(f)(x) = J^{s_1+s_2}(f)(x) = J^{s_2}J^{s_1}(f)(x) \quad (3.2)$$

Proof

$$J^{s_1}J^{s_2}(f)(x) =$$

$$\frac{1}{\Gamma(s_1)\Gamma(s_2)}\int_0^x dy\,(x-y)^{s_1-1}\int_0^y dz\,(y-z)^{s_2-1}f(z) \quad (3.3)$$

we apply the Dirichlet's formula given by Whittaker and Watson [7] [8]

$$\int_0^x dy\,(x-y)^{s_1-1}\int_0^y dz\,(y-z)^{s_2-1}g(y,z)$$

$$= \int_0^x dz\int_x^z dy\,(x-y)^{s_1-1}(y-z)^{s_2-1}g(y,z) \quad (3.4)$$

for $g(y,z) = f(z)$

$$J^{s_1}J^{s_2}(f)(x)$$
$$= \frac{1}{\Gamma(s_1)\Gamma(s_2)}\int_0^x dzf(z)\int_x^z dy\,(x-y)^{s_1-1}(y-z)^{s_2-1}$$

we perform the change of variable

$$u = \frac{y-z}{x-z}$$

$$y = z + u(x-z) \quad dy = (x-z)du$$

$$J^{s_1}J^{s_2}(f)(x)$$
$$= \frac{1}{\Gamma(s_1)\Gamma(s_2)}\int_0^x dzf(z)\,(x-z)^{s_1+s_2}\int_0^1 du\,(1-u)^{s_1-1}u^{s_2-1}$$

we obtain the Euler's beta function

$$B(s_1)(s_2) = \frac{\Gamma(s_1)\Gamma(s_2)}{\Gamma(s_1+s_2)}$$

$$J^{s_1}J^{s_2}(f)(x) = \frac{1}{\Gamma(s_1+s_2)}\int_0^x dzf(z)\,(x-z)^{s_1+s_2-1}$$

$$= J^{s_1+s_2}(f)(x)$$

$J^{s_1}J^{s_2}$ is symmetric in $(s_1, s_2)$, then

$$J^{s_1}J^{s_2}(f)(x) = J^{s_1+s_2}(f)(x) = J^{s_2}J^{s_1}(f)(x)$$

we obtain the semi-group property of $J^s$, which stands true for any function $f$ satisfying the existence of $J^s(f)(x)$ (for any $s \in \mathbb{C}$ with the condition $Re(s) > 0$).

## IV. Definition of $s$-order derivative operator $D^s$ for a complex number $s$ with $Re(s) > 0$

By considering the fact that the operator of one order derivative $D^1$ is a left handside inverse [9] of the operator $J^1$ and extending the relation (1.3), we may give the following definition

### Definition 2

For all function $f$ fulfilling the condition $f(x) = 0$ for $x \leq 0$

and for a complex number $s$ with positive real part, the action of a $s-$order derivative operator $D^s$ on $f$ is defined by the relation

$$D^s(f)(x) = D^k J^{k-s}(f)(x) \qquad (4.1)$$

A sufficient condition which determines the choice of $k$ is that $J^{k-s}(f)(x)$ is well defined.
According to the Theorem 1, we have to choose $Re(k-s) > 0$.
Practically we choose $k \in \mathbb{N}$ and $k > Re(s)$. For instance $k = [Re(s)] + 1$, $[Re(s)]$ is the entire part of $Re(s)$, then

$$D^s(f)(x) = D^{[Re(s)]+1} J^{[Re(s)]+1-s}(f)(x)$$

**Theorem 4**

We have the semi-group property

$$D^{s_1} D^{s_2}(f)(x) = D^{s_1+s_2}(f)(x) = D^{s_2} D^{s_1}(f)(x) \qquad (4.2)$$

for any $s_1$ and $s_2$ with $Re(s_1) > 0$ and $Re(s_2) > 0$

*Proof*
The proof may be easily derived from the semi-group property of $J^s$ (Theorem 3)

*Example 4*
Let us consider $f(x) = x^{a+ib}$
According to the relation (2.4), we have

$$J^{k-s}(x^{a+ib}) = \frac{\Gamma(a+1+ib)}{\Gamma(k-s+a+ib+1)} x^{k-s+a+ib}$$

$$\begin{aligned} D^s(x^{a+ib}) &= D^k J^{k-s}(f)(x) \\ &= D^k \left[ \frac{\Gamma(a+1+ib)}{\Gamma(k-s+a+ib+1)} x^{k-s+a+ib} \right] \\ &= \frac{\Gamma(a+1+ib)}{\Gamma(k-s+a+ib+1)} D^k[x^{k-s+a+ib}] \\ &= \frac{\Gamma(a+1+ib)}{\Gamma(k-s+a+ib+1)} \\ &\quad \cdot \frac{\Gamma(k-s+a+ib+1)}{\Gamma(k-s+a+ib+1-k)} x^{k-s+a+ib-k} \\ &= \frac{\Gamma(a+1+ib)}{\Gamma(-s+a+ib+1)} x^{-s+a+ib} \qquad (4.3) \end{aligned}$$

The result is independent on $k$

If $s = \alpha + i\beta$

$$D^s(x^{a+ib}) = \frac{\Gamma(a+ib+1)}{\Gamma(-s+a+ib+1)} x^{a+ib-s}$$

$$= \frac{\Gamma(a+ib+1)}{\Gamma(-\alpha-i\beta+a+ib+1)} x^{a-\alpha+i(b-\beta)}$$

$$= \frac{\Gamma(\alpha+1+i\beta)}{\Gamma(a-\alpha+1+ib-i\beta)} x^{a-\alpha}$$

$$\cdot [cos((b-\beta)lnx) + isin((b-\beta)lnx)] \qquad (4.4)$$

From the relation (4.3), it can be easily shown that

$$D^{s_1} D^{s_2}(x^{a+ib}) = D^{s_1+s_2}(x^{a+ib}) = D^{s_2} D^{s_1}(x^{a+ib}) \qquad (4.5)$$

Comparing the expressions (4.3) and (2.4), we may write the relations

$$J^{-s}(x^{a+ib}) = D^s(x^{a+ib}) \qquad (4.6a)$$

$$D^{-s}(x^{a+ib}) = J^s(x^{a+ib}) \qquad (4.6b)$$

These relations suggest then the extension of the definition of the operators $J^s$ and $D^s$ for any complex number $s$.

## V. DEFINITION OF $s$-ORDER INTEGRAL OPERATOR $J^s$ AND $s$-ORDER DERIVATIVE OPERATOR $D^s$ FOR ANY COMPLEX NUMBER

The results (4.6) suggest the following definition

- If $Re(s) > 0$

$$J^s(f)(x) = \frac{1}{\Gamma(s)} \int_0^x (x-y)^{s-1} f(y) dy \qquad (5.1a)$$

$$= \frac{x^s}{\Gamma(s)} \int_0^1 (1-u)^{s-1} f(ux) du$$

- If $Re(s) \leq 0$

$$J^s(f)(x) = D^k J^{k+s}(f)(x) \qquad (5.1b)$$

$$= D^k \left[ \frac{1}{\Gamma(s)} \int_0^x (x-y)^{k+s-1} f(y) dy \right]$$

and

$$D^s(f)(x) = J^{-s}(f)(x) \qquad (5.2)$$

The constant $k$ in the relation $(5.1b)$ is to be chosen such $J^{k+s}(f)(x)$ is well defined. According to the Theorem 1, we have to choose $Re(k+s) > 0$. Practically we choose $k \in \mathbb{N}$ and $k > Re(-s)$. For instance, one can choose $k = [Re(-s)] + 1$, $[Re(-s)]$ is the entire part of $Re(-s)$.

## VI. EXTENSION OF THE DEFINITION DOMAIN. DEFINITION OF $s$-ORDER INTEGRAL OPERATOR $J^s_{x_0}$ AND $s$-ORDER DERIVATIVE OPERATOR $D^s_{x_0}$ FOR A COMPLEX NUMBER $s$ WITH $Re(s) > 0$

For a given real number $x_0$, the definition of $s$-order integral may be extended for more general function with the introduction of the operator $J^s_{x_0}$ defined in the relation (6.1). The action of the operator $J^s_{x_0}$ on a function $f$ is defined firstly, for complex number $s$ with positive real part. Then, using the fact that the operator of one order derivative, noted $D^1$, is the left hand side inverse of the operator $J_{x_0}^s$, an $s$-order derivative operator, noted $D_{x_0}^s$ is also defined for all complex number $s$ with positive real part. Finally, assuming the relation $J_{x_0}^s = D_{x_0}^{-s}$ the definition of the $s$-order integral and $s$-order derivative is extended for any complex number $s$ with positive and negative real parts.

**Theorem 5**

For a derivable and integrable function $f$ defined on the interval $I = ]x_0, +\infty[$, $x_0 \in \mathbb{R}$, the integral $J_{x_0}^s(f)(x)$ defined by the relation

$$J_{x_0}^s(f)(x) = \frac{1}{\Gamma(s)} \int_{x_0}^{x} (x-y)^{s-1} f(y) dy \qquad (6.1)$$

is absolutely convergent for $s \in \mathbb{C}$ with $Re(s) > 0$.

*Proof*

Let us take $s = \alpha + i\beta$; It is to be noted that $x_0 \leq y \leq x$ so $(x-y)^\alpha \geq 0$, then

$$\left| \frac{(x-y)^{s-1}}{\Gamma(s)} f(y) \right| = \frac{|x-y|^{\alpha-1}}{|\Gamma(s)|} |f(y)|$$

$$\Rightarrow \int_{x_0}^{x} \left| \frac{(x-y)^{\alpha-1}}{\Gamma(s)} f(y) \right| dy = \int_{x_0}^{x} \frac{(x-y)^{\alpha-1}}{|\Gamma(s)|} |f(y)| dy$$

Using an integration by parts, we have for $Re(s) = \alpha > 0$

$$\int_{x_0}^{x} \frac{(x-y)^{\alpha-1}}{|\Gamma(s)|} |f(y)| dy$$

$$= -\frac{1}{\alpha |\Gamma(s)|} \left\{ \lim_{y \to x^-} (x-y)^\alpha |f(y)| - \lim_{y \to x_0^+} (x-y)^\alpha |f(y)| - \int_{x_0}^{x} (x-y)^\alpha |f(y)|' dy \right\}$$

$$\lim_{y \to x^-} (x-y)^\alpha |f(y)| = 0$$

$$\lim_{y \to x_0^+} (x-y)^\alpha |f(y)| = (x-x_0)^\alpha \lim_{y \to x_0^+} |f(y)|$$

So

$$\int_{x_0}^{x} \frac{(x-y)^{\alpha-1}}{|\Gamma(s)|} |f(y)| dy = \frac{(x-x_0)^\alpha}{\alpha |\Gamma(s)|} \lim_{y \to x_0^+} |f(y)|$$

$$+ \int_{x_0}^{x} \frac{(x-y)^\alpha}{\alpha |\Gamma(s)|} |f(y)|' dy$$

For $x \in I = ]x_0, +\infty[$ and for $\alpha > 0$, we have

$$x_0 \leq y \leq x \;\; \Rightarrow 0 \leq x - y \leq x \Rightarrow 0 \leq (x-y)^\alpha \leq x^\alpha$$

$$\Rightarrow 0 \leq \frac{(x-y)^\alpha}{\alpha |\Gamma(s)|} \leq \frac{(x-x_0)^\alpha}{\alpha |\Gamma(s)|}$$

$$\Rightarrow \int_{x_0}^{x} \frac{(x-y)^\alpha}{\alpha |\Gamma(s)|} |f(y)|' dy \leq \frac{(x-x_0)^\alpha}{\alpha |\Gamma(s)|} \int_{x_0}^{x} |f(y)|' dy$$

$$\int_{x_0}^{x} \frac{(x-y)^\alpha}{\alpha |\Gamma(s)|} |f(y)|' dy \leq \frac{(x-x_0)^\alpha}{\alpha |\Gamma(s)|} [|f(x)| - \lim_{y \to x_0^+} |f(y)|]$$

then

$$\int_{x_0}^{x} \left| \frac{(x-y)^{s-1}}{\Gamma(s)} f(y) \right| dy = \frac{(x-x_0)^\alpha}{\alpha |\Gamma(s)|} \lim_{y \to x_0^+} |f(y)|$$

$$+ \int_{x_0}^{x} \frac{(x-y)^\alpha}{\alpha |\Gamma(s)|} |f(y)|' dy \leq \frac{(x-x_0)^\alpha}{\alpha |\Gamma(s)|} \lim_{y \to x_0^+} |f(y)| +$$

$$\frac{(x-x_0)^\alpha}{\alpha |\Gamma(s)|} [|f(x)| - \lim_{y \to x_0^+} |f(y)|]$$

Then for $\alpha = Re(s) > 0$

$$\int_{x_0}^{x} \left| \frac{(x-y)^{s-1}}{\Gamma(s)} f(y) \right| dy = \int_{x_0}^{x} \frac{(x-y)^\alpha}{\alpha |\Gamma(s)|} |f(y)|' dy$$

$$\leq \frac{(x-x_0)^\alpha}{\alpha |\Gamma(s)|} |f(x)|$$

The integral $J_{x_0}^s(f)(x)$ is absolutely convergent for $\alpha = Re(s) > 0$.

**Theorem 6**

Let $D^1$ be the operator of one –order derivative

$$D^1 = \frac{d}{dx}$$

For any complex number $s$ with $Re(s) > 0$, we have

$$D^1\left(J_{x_0}^s\right)(f)(x) = J_{x_0}^{s-1}(f)(x) \qquad (6.2)$$

*Proof*

$$D^1(J_{x_0}^s)(f)(x) = D^1\left[\frac{1}{\Gamma(s)}\int_{x_0}^x (x-y)^{s-1} f(y) dy\right]$$

$$= \frac{(s-1)}{\Gamma(s)}\int_{x_0}^x (x-y)^{s-2} f(y) dy$$

$$= \frac{1}{\Gamma(s-1)}\int_{x_0}^x (x-y)^{s-2} f(y) dy = J_{x_0}^{s-1}(f)(x)$$

From this theorem, it can be easily shown that

- the operator $D^1$ is a left hand side inverse of the operator $J_{x_0}^1$

- for $n, k \in \mathbb{N}, n \geq k$

$$D^1(J_{x_0}^s)(f)(x) = J_{x_0}^{s-1}(f)(x) \qquad (6.3)$$

These results justify the introduction of the following definition

*Definition 4*

Let $f$ be a derivable and integrable function defined on the interval $I = ]x_0, +\infty[$. For a complex number $s$ with $Re(s) > 0$, we define an $s$-order integral of the function $f$, noted $J_{x_0}^s$, by the relation

$$J_{x_0}^s(f)(x) = \frac{1}{\Gamma(s)}\int_{x_0}^x (x-y)^{s-1} f(y) dy \qquad (6.4)$$

According to the theorem 5, the integral $J_{x_0}^s(f)(x)$ is well defined.

*Example 6.1*

Let us take $x_0 = 0$ and $f(x) = x^p$ with $p \in \mathbb{C}$. From the relation (2.4), we have

$$J_0^s(f)(x) = \frac{1}{\Gamma(s)}\int_0^x (x-y)^{s-1} f(y) dy$$

$$= \frac{\Gamma(s)\Gamma(p+1)}{\Gamma(s+p+1)}\frac{x^{s+p}}{\Gamma(s)} = \frac{\Gamma(p+1)}{\Gamma(s+p+1)} x^{s+p} \qquad (6.5)$$

*Example 6.2*

Let us take $x_0 \to -\infty$ and $f(x) = e^x$. Then

$$J_{-\infty}^s(e^x) = \frac{1}{\Gamma(s)}\int_{-\infty}^x (x-y)^{s-1} e^y dy$$

For instance, for $s = 1$, we have

$$J_{-\infty}^1(e^x) = \int_{-\infty}^x e^y dy = e^x$$

Using mathematical induction, it may be proven easily that

$$J_{-\infty}^n(e^x) = e^x \text{ for any } n \in \mathbb{N}$$

*Definition 5*

Let $f$ be a derivable and integrable function defined on the interval $I = ]x_0, +\infty[$. We may define an $s$-order derivative operator $D_{x_0}^s$ for a complex number $s$ by using the fact that the operator of one order derivative $D^1$ is a left hand side inverse of the operator $J_{x_0}^1$

$$D_{x_0}^s(f)(x) = D^k J_{x_0}^{k-s}(f)(x) \qquad (6.6)$$

A sufficient condition which determines the choice of $k$ is that $J_{x_0}^{k-s}(f)(x)$ is well defined. According to the relation (5.2), we have to choose $Re(k-s) > 0$. Practically, we choose $k \in \mathbb{N}$ and $k > Re(s)$, for instance $k = [Re(s)] + 1$, in which $[Re(s)]$ is the entire part of $[Re(s)]$.

$$D_{x_0}^s(f)(x) = D^{[Re(s)]+1} J^{[Re(s)]+1-s}(f)(x) \qquad (6.7)$$

*Example 6.3*

Let us study the example of $f(x) = x^p$ and $x_0 = 0$, According to (6.6), we have

$$J_0^{k-s}(x^p) = \frac{\Gamma(k-s)\Gamma(p+1)}{\Gamma(k-s+p+1)}\frac{x^{k-s+p}}{\Gamma(k-s)}$$

$$= \frac{\Gamma(p+1)}{\Gamma(k-s+p+1)} x^{k-s+p}$$

So

$$D_0^s(x^p) = D^k\left[\frac{\Gamma(p+1)}{\Gamma(k-s+p+1)} x^{k-s+p}\right]$$

$$= \frac{\Gamma(p+1)}{\Gamma(-s+p+1)} x^{-s+p}$$

## VII. CONCLUSIONS

According to the results that we have obtained, particularly the relations (5.1a), (5.1b), (6.4) and (6.5), we may conclude that the approach that we have considered allows the extension of the definitions of $s$-order integrals and $s$-order derivatives for the case of complex order $s$. Generally speaking, we obtain the same formal expressions as in the case of $s \in \mathbb{R}$ [4]. We have considered the case of functions which fulfil more general condition (paragraph VI) too.

Deep analysis of the above results can be done to obtain a better understanding about the fundamental and physical meaning of integrals and derivatives for the case of complex order.